\documentclass[letterpaper, 10 pt, journal, twoside]{IEEEtran}
\usepackage{amsmath,amsfonts}
\usepackage{algorithmic}
\usepackage{algorithm}
\usepackage{array}
\usepackage{subfigure}
\usepackage{textcomp}
\usepackage{stfloats}
\usepackage{url}
\usepackage{verbatim}
\usepackage{graphicx}
\usepackage{cite}
\usepackage{amsthm}
\usepackage{multirow}

\theoremstyle{definition}

\begin{document}

\title{UAV Routing for Enhancing the Performance of a Classifier-in-the-loop}

\author{Deepak Prakash Kumar$^{1}$, Pranav Rajbhandari$^{2}$, Loy McGuire$^{3}$, Swaroop Darbha$^{4}$, Donald Sofge$^{5}$
\thanks{$^{1}$ Deepak Prakash Kumar is a graduate student in the Department of Mechanical Engineering, Texas A\&M University (e-mail: {\tt\footnotesize deepakprakash1997@gmail.com}).}%
\thanks{$^{2}$ Pranav Rajbhandari is an undergraduate student in the Department of Computer Science, Carnegie Mellon University (e-mail: {\tt\footnotesize prajbhan@andrew.cmu.edu})}%
\thanks{$^{3}$ Loy McGuire is a graduate student in the Department of Aerospace Engineering, University of Maryland (e-mail: {\tt\footnotesize loy.mcguire@nrl.navy.mil})}%
\thanks{$^{4}$ Swaroop Darbha is a professor in the Department of Mechanical Engineering, Texas A\&M University (e-mail:
{\tt\footnotesize dswaroop@tamu.edu})}%
\thanks{$^{5}$ Donald Sofge is the section head of Distributed Autonomous Systems, Navy Center for Applied Research in Artificial Intelligence, US Naval Research Laboratory, Washington (e-mail:
{\tt\footnotesize donald.sofge@nrl.navy.mil})}%
}


\markboth{IEEE TASE. Submission Version.}%
{Kumar \MakeLowercase{\textit{et al.}}: UAV Routing for Enhancing Performance of Classifier-in-the-loop}


\maketitle

\begin{abstract}
Some human-machine systems are designed so that machines (robots) gather and deliver data to remotely located operators (humans) through an interface in order to aid them in classification. The performance of a human as a (binary) classifier-in-the-loop is characterized by probabilities of correctly classifying objects of type $T$ and $F$. These two probabilities depend on the dwell time, $d$, spent collecting information at a point of interest (POI or interchangeably, target). The information gain associated with collecting information at a target is then a function of dwell time $d$  and discounted by the revisit time, $R$, i.e., the duration between consecutive revisits to the same target. The objective of the problem of routing for classification is to optimally route the vehicles and determine the optimal dwell time at each target so as to maximize the total discounted information gain while visiting every target at least once. In this paper, we make a simplifying assumption that the information gain is discounted exponentially by the revisit time; this assumption enables one to decouple the problem of routing with the problem of determining optimal dwell time at each target for a single vehicle problem. For the multi-vehicle problem, we provide a fast heuristic to obtain the allocation of targets to each vehicle and the corresponding dwell time.
\end{abstract}

\def\abstractname{Note to Practitioners}
\begin{abstract}
The motivation for the considered problem is primarily for surveillance applications, wherein certain locations/points of interest (POIs) need to be monitored for activity. An operator/classifier will identify such points and provide them to an interface. Through the proposed model and algorithm in this study, unmanned aerial vehicles (UAVs) will be provided with the route to be taken and targets to be covered to collect information on these points of interest. Once such information reaches the operator, the operator can decide whether each POI is a target, not a target, or whether additional information is required.
\end{abstract}


\section{Introduction}

\IEEEPARstart{T}{he} design of human-machine systems requires a careful partition of the tasks to be performed by the human-in-the-loop and the machines, and the design of an associated human-machine interface. In this paper, we consider a human-machine system, where the human serves as a classifier-in-the-loop based on the information delivered to the human operator by the machines (vehicles) through an interface.
The interface takes as input $n$ Points of Interest (POIs/targets) nominated by the human operator and computes the order in which they must be visited, the time to be spent at each POI along with the waypoints for the vehicles to follow. The vehicles persistently monitor the targets by visiting them and dwelling at a POI while collecting information that is transmitted to the remotely located operator. Based on the information supplied by the vehicles, the primary task of the human is to assess/classify events happening at the $n$ specified targets and classify them as $T$ (true target) or $F$ (false target). The probability of an operator correctly classifying events at a target depends on the dwell time of the vehicle at that location and the revisit time, $R$, i.e., the time duration between consecutive revisits to the same location. A discounted information gain reflects the tradeoff between the information gained as a function of dwell time $d$ and the revisit time, $R$. The problem considered in this paper is to determine the optimal sequence in which targets must be visited and the dwell time at each POI that helps the operator maximize the sum of all discounted information gained from targets, while ensuring that each of them is visited. While this problem is solved by the interface; in a centralized manner, distributed heuristics can be developed for this problem; however, it is not considered in this paper. Subsequent determination of waypoints for vehicles follows immediately and is communicated by the interface to the vehicles.

\begin{figure}[htb!]
    \centering
    \includegraphics[width = 0.85\linewidth]{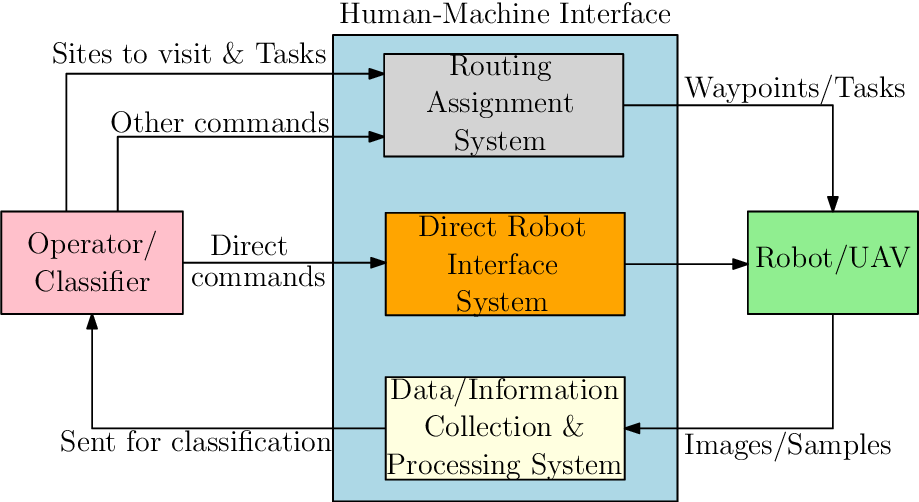}
    \caption{Human-Machine System with human as a classifier-in-the-loop}
    \label{fig: human_machine_interface}
\end{figure}

The novel contributions of this paper are as follows:
\begin{itemize}
    \item A novel model for accounting the dependence of UAV's dwell time on information gain, and choosing an objective for routing based on the model to enhance the operator's performance as a classifier-in-the-loop.
    \item The objective of routing is to maximize information gain, a nonlinear function of the continuous decision variables (dwell times) and discrete decision variables (routing choice); corresponding heuristics to solve this difficult routing problem and associated parametric numerical study are other novel contributions. 
\end{itemize}

\section{Literature Review}
\label{Literature Review}
Consider a set of $n$ POIs/targets to be visited by $m$ vehicles to gain information about each target. Let $G = (V, E)$ be an undirected graph where $V$ is the set of targets to be visited and $E$ is the set of edges with weights corresponding to known non-negative traveling costs between the targets. These costs need not necessarily be the Euclidean distances between targets and may depend on other factors such as environmental considerations. It will be assumed the cost of traveling from $t$ to $u$ is the same as the cost to travel from $u$ to $t$. Each vehicle departs and returns to a depot $D_j$ for $j = 1, 2, \cdots, m$. It is desired for each target to be visited by exactly one vehicle. Any vehicle visiting a POI collects information and transmits it to a remotely located operator in order to aid its classification. The means of gaining this information by the vehicle is not considered in this paper. Each target is to be classified as either $T$ (target) or $F$ (not a target) with the exception of the depots which will have no classification. It will be assumed there is no \textit{a priori} information about the targets. It will also be assumed the probability of correctly classifying the POI is the same whether the target has a true classification of $T$ or $F$. That is, it is equally difficult to classify the target, regardless of what it really is. The objective is to construct tours for the vehicles such that each target is exactly visited by one target, and the total information gained is maximized.

The topic of routing unmanned vehicles for autonomously collecting data has received significant attention in the literature \cite{Rathinam_IEEE_TASE2007, hari_IEEETOR_2021_550,hari2022bounds, hari2022cooperative, MALIK_ORL_2007_747, manyam2015computation, manyam2017tightly, Manyam_IEEERAL_2022, Oberlin_IEEERAM_2010, oberlin2009transformation, obermeyer_JGCD_2012sampling, bhadoriya2022assisted, yadlapalli_OL_20123}. The objective of routing depends on the mission, and the nature of the optimal solution depends further on the operational, motion, coordination and communication constraints. While the previous references addressed routing problems, none of them deal with routing vehicles to enhance the performance of an operator as a classifier-in-the-loop. Typically, the performance of an operator in the loop is specified by a confusion matrix, which specifies two conditional probability distributions -- given that the event/object is of type $Z \in \{T, F\}$, the conditional probability distribution specifies the discrete probability distribution of correctly classifying the event/object\cite{kalyanam2016optimal}. These probability distributions depend on controllable operational parameters such as altitude, pose of the vehicle relative to the object and the time spent imaging the object/event, etc. The idea of vehicles enhancing the classification performance then rests on controlling or choosing these parameters to ensure that the conditional distributions are separated as far as possible, i.e, the mutual information gain is maximized. For persistent monitoring of the targets, it is also necessary to minimize the time, $R$, between successive revisits to the targets. This is to ensure that discounting the mutual gain by a function of the revisit time is reasonable and the corresponding objective for optimization is to maximize the discounted mutual information gain over all targets through optimal routing and determination of optimal dwell time.

Information gain\cite{information_theory_book} has been used in path planning in robotic applications that are distinct from what is considered in this paper. 
Lee \textit{et al.} \cite{lee2010enhanced} developed an enhanced ant colony optimization for the capacitated vehicle routing problem by using information gain to ameliorate the search performance when the good initial solution was provided by simulated annealing algorithm. Toit and Burdick \cite{du2011robot} used the information gain theory in developing a partially closed-loop receding horizon control algorithm to solve the stochastic dynamic programming problem associated with dynamic uncertain environments robot motion planning. Kaufman \textit{et al.} \cite{kaufman2016autonomous} presented a novel, accurate and computationally-efficient approach to predict map information gain for autonomous exploration where the robot motion is governed by a policy that maximises the map information gain within its set of pose candidates. Zaenker \textit{et al.} \cite{zaenker2023graph} proposed a novel view motion planner for pepper plant monitoring while minimizing occlusions (a significant challenge in monitoring of large and complex structures), that builds a graph network of viable view poses and trajectories between nearby poses which is then searched by planner for graphs for view sequences with highest information gain. Paull \textit{et al.} \cite{paull2010information} used information gain approach in the objective function of sidescan sonars (SSS) and for complete coverage and reactive path planning of an autonomous underwater vehicle. Mostofi \cite{mostofi2009decentralized} proposed a communication-aware motion-planning strategy for unmanned autonomous vehicles, where each node considers the information gained through both its sensing and communication when deciding on its next move. They showed how each node can predict the information gained through its communications, by online learning of link quality measures and combining it with the information gained through its local sensing in order to assess the overall information gain.\par

Information gain finds its place in machine learning literature where it is being used for diverse feature ranking and feature selection techniques in order to discard irrelevant or redundant features from a given feature vector, thus reducing dimensionality of the feature space. Novakovic \cite{novakovic2009using} applied Information Gain for the classification of sonar targets where the Information Gain evaluation helped in increasing computational efficiency while improving classification accuracy by doing feature selection.\par
Information-theoretic methods have been used for computing heuristics for path-planning methods in autonomous robotic exploration where mutual information is calculated between the sensor's measurements and the explored map. Deng \textit{et al.} \cite{deng2020robotic} proposed a novel algorithm for the optimizing exploration paths of a robot to cover unknown 2D areas by creating a gradient-based path optimization method that tries to improve path's smoothness and information gain of uniformly sampled view-points along the path simultaneously. Julian \textit{et al.} \cite{julian2014mutual} proved that any controller tasked to maximise a mutual information reward function is eventually attracted to unexplored space which is derived from the geometric dependencies of the occupancy grid mapping algorithm and the monotonic properties of mutual information. Bai \textit{et al.} \cite{bai2016information} proposed a novel approach to predict mutual information using Bayesian optimisation for the purpose of exploring a priori unknown environments and producing a comprehensive occupancy map. They showed that information-based method provides not only computational efficiency and rapid map entropy reduction, but also robustness in comparison with competing approaches. Amigoni \& Caglioti \cite{amigoni2010information} presented a mapping system that builds geometric point-based maps of environments employing an information-based exploration strategy that determines the best observation positions by blending together expected gathered information (that is measured according to the expected  a posteriori uncertainty of the map) and cost of reaching observation positions. Basilico \& Amigoni \cite{basilico2011exploration} further extended this information-based exploration strategy for rescue and surveillance applications. In \cite{Alipour2011OnDO}, the problem of routing a mobile agents for data aggregation in sensor networks is considered. Here, the main issue is the tradeoff between increasing information gain and power consumption among the source nodes that must be visited by the mobile agent and is accounted for in the cost of the edges.

The problem of routing vehicles for aiding an operator-in-the-loop for classification was first proposed by Montez \cite{Montez_2020_SPIE, SuveerMSThesis2023}; however, the paper \cite{Montez_2020_SPIE} does not exploit the exponential discounting nature of mutual information gain to decouple the mixed-integer nonlinear program into a discrete optimization problem and a continuous optimization problem. The form for $P_i$ considered in \cite{SuveerMSThesis2023} does not possess the desired structure for information gain that we seek in this article and is different from \cite{SuveerMSThesis2023} in that respect.  This structure is exploited in this paper and an exact algorithm for a single vehicle routing is presented in this paper. In addition, extensions to the multiple vehicle case is presented with some heuristics along with the corroborating computational results.

\section{Mathematical Formulation}

\subsection{Quantifying the Information Gained}

Suppose a vehicle visits the $i\textsuperscript{th}$ target. Denote the set of classification choices as $C = \{T, \, F\}$. Each POI has a correct classification $X \in C$. The operator assigns a classification of $Z \in C$ to the $i\textsuperscript{th}$ POI after the visit. Let $s_i$ represent the variables affecting the observation. Denote the conditional probabilities of correctly classifying $i$ as $T$ or $F$ given the variables $s_i$ as
\begin{align*}
    P_t(s_i) &= P(Z = T \, | \, X = T, \, s_i) \, \, \mbox{and} \\
    P_f(s_i) &= P(Z = F \, | \, X = F, \, s_i), 
\end{align*}
respectively. Associated with the classification, we have a confusion matrix given below:
\begin{table}[htb!]
\caption{Confusion Matrix}
\begin{center}
\begin{tabular}{||c||c|c||}
\hline \hline  
& $Z=T$& $Z=F$\\
\hline \hline 
$X=T$& $P_t(s_i)$ & $1-P_t(s_i)$ \\
\hline 
$X=F$ & $1-P_f(s_i)$& $P_f(s_i)$\\
\hline \hline 
\end{tabular}
\end{center}
\end{table}
The two rows of the confusion matrix indicate the probability distribution of classification conditioned on the POI being of type $T$ and $F$ respectively, and depend on the two parameters $P_t(s_i)$ and $P_f(s_i)$. 
Ideally, one would want to separate the two conditional probability distributions as much as possible with controllable variables which affect the observation, such as dwelling more time at a POI or observing the POI at a lower altitude or from a better perspective. The information gained by visiting each POI can be quantified using the Kullback-Leibler divergence (also referred to as the mutual information or information gain) as the distance between the two conditional probability distributions. The mutual information for $i \in V$ between the two classification variables $X$ and $Z$ will be denoted as $I_i(X, Z)$. The mutual information is defined to be
\begin{align} \label{Mutual Information Definition}
    I_i = \sum_{x,z \in C} P(X = x, Z = z) \log{\frac{P(X = x, Z = z)}{P(X = x) P(Z = z)}}.
\end{align}
Denote the \textit{a priori} probability a POI is a target, $P(X = T),$ as $p$. 
It will be assumed the \textit{a priori} probability a POI is a target is 0.5. That is, there is effectively no known information about the targets before sending out the vehicle to investigate, and so each POI is equally likely to be either a target or not a target. Additionally, it will be assumed it is equally difficult to correctly classify the $i\textsuperscript{th}$ POI, as a target or not a target. That is, $P_t(s) = P_f(s) = P_i(s)$ for any set of variables $s_i$ that affect observation. Then, \eqref{Mutual Information Definition} reduces to
\begin{align} \label{Reduced Mutual Information}
\begin{split}
    I_i &= P_i(s_i) \log P_i(s_i) + (1 - P_i(s_i)) \log(1 - P_i(s_i)) \\
    & \quad\, + \log 2.
\end{split}
\end{align}
\begin{figure}[htb!]
    \centering
    \includegraphics[width = 0.6\linewidth]{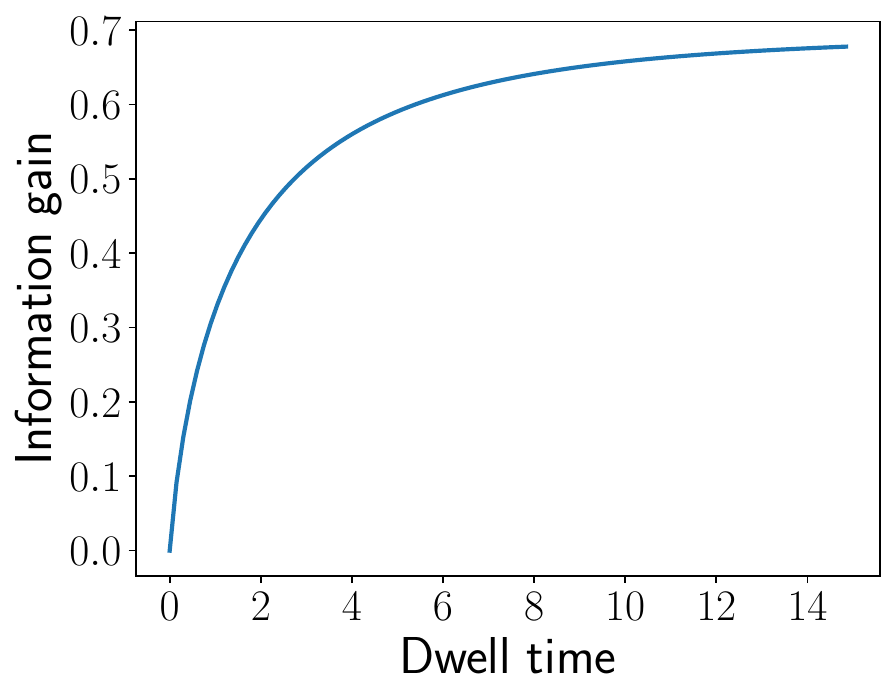}
    \caption{Information Gain vs  Dwell Time at a POI ($\tau = 0.5$)}
    \label{fig:infogain1}
\end{figure}
If $P_i(s_i) = P_i(d_i)$, a function of the dwell time, $d_i$, at the $i\textsuperscript{th}$ POI, then one can express mutual information gain $I_i$ at the $i\textsuperscript{th}$ POI as an explicit function of the dwell time $d_i$. At this point, we observe the following properties of the mutual information gain function (see Fig.~\ref{fig:infogain1}, where $\tau$ is a parameter to be defined later):
\begin{itemize}
    \item The function $I_i(d_i)$ is monotonically increasing with $d_i$; essentially, the information gain increases with the time spent by a vehicle at the $i\textsuperscript{th}$ target.  Hence, $\frac{\partial I_i}{\partial d_i} \ge 0$. 
    \item Law of diminishing returns applies to the information gain, i.e., the marginal increase in information gain decreases with the dwell time. Hence, $\frac{\partial^2 I_i}{\partial d_i^2} \le 0$. 
    \item Information gained is always non-negative, i.e., $I_i(d_i) \ge 0$. 
\end{itemize}

A consequence of these properties is that $I_i(d_i)$ is log-concave as $I_i(d_i)\frac{\partial^2 I_i}{\partial d_i^2} - (\frac{\partial I_i}{\partial d_i} )^2 \le 0 $. It is also true that 
$J_0 = \sum_{i=1}^n I_i(d_i)$
is log-concave since 
$$\left[\sum_{i=1}^n I_{i}(d_i)\right]\left[ \sum_{j=1}^n \frac{\partial^2 I_j}{\partial d_j^2}\right]- \left[\sum_{i=1}^n \frac{\partial I_i}{\partial d_i}\right]^2 \le 0. $$ 
A consequence of this observation is that one may employ gradient ascent to $log (J_0 (d_1, \ldots, d_n))$ to arrive at the optimum.
In this paper, we model $P_i(s_i)$ as
\begin{equation} \label{P_i definition}
    P_i(s_i) = P_i(d_i) = 1- \frac{1}{2} e^{-\sqrt{d_i/\tau_i}} \, ,
\end{equation}
where $\tau_i$ is a positive constant that represents the sensitivity to the time spent at the  $i\textsuperscript{th}$ target. 
This form of $P_i(s_i)$ is consistent with the observations listed above for information gain. 
Correspondingly, the information gain may be expressed as solely a function of $d_i$ as
\begin{align*}
    I_i(d_i) &= \left(1-\frac{1}{2}e^{-\sqrt{d_i/\tau_i}}\right) \log\left(1-\frac{1}{2}e^{-\sqrt{d_i/\tau_i}}\right) \\
    & \quad\, - \frac{1}{2}e^{-\sqrt{d_i/\tau_i}}\left(\log 2 +\sqrt{d_i/\tau_i}\right) + \log 2.
\end{align*}
A sample plot of information gain corresponding to $\tau_i = 0.5$ is given in Fig.~\ref{fig:infogain1}. 

Since we want to incentivize the vehicles to visit all targets, we discount the information gain by the revisit time, $R_i$, for the $i\textsuperscript{th}$ target as follows:
$$ \psi_i(d_i, R_i) = e^{-\alpha R_i} I_i(d_i),$$
where $\alpha>0$ is a positive constant, $R_i$ is the time duration between successive revisits to the $i\textsuperscript{th}$ target. 

The objective of the optimization problem considered in this paper is to maximize
\begin{align} \label{eq: objective_function_single_vehicle}
    J_s (d_1, d_2, \ldots, d_n) = \sum_{i=1}^n \psi_i(d_i, R_i),
\end{align}
through the choice of a route for the vehicles and the dwell time at each POI, while ensuring that each target is visited. We note here that since $I_i$ is log-concave, $J_s$ is log-concave as well.

\subsection{Single Vehicle Case}
In the case of a single vehicle, $R_i$ is the same for every target (say, it is $R$) if every other target is visited exactly once between successive revisits; moreover $R = T + \sum_{i=1}^n d_i$, where $T$ is time taken to tour the $n$ targets. If triangle inequality holds, this is true even if one may allow the same target to be visited multiple times between consecutive revisits to another target \cite{hari_IEEETOR_2021_550}. Note that $T \ge TSP^*$, where $TSP^*$ is the minimum time taken to visit the $n$ targets before returning to the starting location. A consequence is the following:
\begin{align*}
   e^{-\alpha R} &\leq e^{-\alpha TSP^*} e^{-\alpha \sum_{i=1}^n d_i}, \\
   \implies J &\leq \sum_{i=1}^n e^{-\alpha TSP^*} e^{ -\alpha \sum_{i=1}^nd_i} I_i(d_i), \\
   &\leq e^{-\alpha TSP^*}  \max_{d_1, \ldots, d_n} e^{ -\alpha\sum_{i=1}^nd_i} \sum_{i=1}^n  I_i(d_i).
\end{align*}
If $J^*$ is the optimum, clearly, it is achieved by minimizing $T$, and maximizing the log-concave function on the right hand side of the above inequality. In other words, the problem of optimal routing and the determination of optimal dwell time at each targets is now decoupled, and can be solved to optimality.
\subsection{Multiple Vehicle Case}
An additional complication arises in the multiple vehicle case -- that of partitioning and assigning the targets to be visited by each vehicle. If there are $m > 1$ vehicles, let the targets be partitioned into $m$ {\it disjoint} sets, namely ${\mathcal P}_1, \ldots, {\mathcal P}_m$, so that the $i\textsuperscript{th}$ vehicle is tasked with visiting the POIs in ${\mathcal P}_i$. Let $R_i$ be the revisit time associated with targets assigned to $i\textsuperscript{th}$ vehicle, and the associated tour cost for persistent monitoring per cycle be $TSP^*({\mathcal P}_i)$. Associated with the $i\textsuperscript{th}$ vehicle, the discounted information gained is given by 
\begin{align} \label{eq: objective_function_one_vehicle_given_partition}
\begin{split}
    J_i (\mathcal{P}_i) &= \max_{d_j, j \in {\mathcal P}_i} e^{-\alpha R_i} \sum_{j \in {\mathcal P}_i} I_j(d_j) \\
    &= e^{-\alpha TSP^*({\mathcal P}_i)} \max_{j \in {\mathcal P}_i} e^{-\alpha \left(\sum_{j\in{\mathcal P}_i}d_j\right)} \sum_{j \in {\mathcal P}_i} I_j(d_j).
\end{split}
\end{align}
Corresponding to the partitions ${\mathcal P}_1, {\mathcal P}_2, \cdots, {\mathcal P}_m,$ the discounted information gain is given by
\begin{align} \label{eq: objective_function_given_partition}
    J (\mathcal{P}_1, \mathcal{P}_2, \cdots, \mathcal{P}_m) = \sum_{i = 1}^{m} J_i (\mathcal{P}_i),
\end{align}
where $J_i (\mathcal{P}_i)$ is given in \eqref{eq: objective_function_one_vehicle_given_partition}. Correspondingly, the objective is to maximize the discounted information gain over all possible partitions, sequences of visiting POIs by every vehicle, and the dwell time at each target:
\begin{align} \label{eq: objective_function_multiple_vehicles}
    J = \max_{{\mathcal P}_i, \; 1\le i\le m} J (\mathcal{P}_1, \mathcal{P}_2, \cdots, \mathcal{P}_m).
\end{align}
Since maximizing over partitions is a difficult combinatorial problem, we provide heuristics for the outer layer of optimization in the above optimization problem and use the single vehicle algorithm for the inner layer of optimization.

\section{Heuristic for Multi-Vehicle Case} \label{sect: heuristic}

From the previous section, it can be observed that the single-vehicle case can be solved to optimality. For the multi-vehicle case, the problem is a coupled routing problem and continuous optimization problem stemming from dwell-time optimization. Hence, in this paper, a heuristic is proposed to obtain high-quality solutions for the multi-vehicle case. To this end, a heuristic inspired by the MD algorithm \cite{MD_algorithm}, a heuristic that yields high-quality solutions for a min-max multi-vehicle multi-depot problem, is discussed. The intuition behind using a similar heuristic structure for a min-max problem for the proposed problem is as follows:
\begin{itemize}
    \item From the objective function, a vehicle having a high tour cost will have a low objective value due to the high penalty incurred due to the exponential term (the discounting term).
    \item Further, a vehicle visiting very few targets will have a low objective value due to very minimal information gain.
\end{itemize}
Hence, similar to the min-max problem, it is desired for vehicle tours to be generated such that the vehicles are load-balanced.

The heuristic is structured through three steps: (i) generation of an initial feasible solution, (ii) local search, and (iii) perturbation of solution. In the following subsections, the three steps in the heuristic will be expanded. For this purpose, the notation used for the graph will be discussed first. Let $T$ denote the set of targets to be covered by $m$ vehicles in a graph $G$. Let the $i\textsuperscript{th}$ vehicle start at depot $D_i$ for $i \in \{1, 2, \cdots, m \}$. Hence, the vertices of graph $G$ are $T \bigcup \{D_1, D_2, \cdots, D_m\}$. We note here that the depots need not be necessarily distinct. The set of edges $E$ in graph $G$ are assumed to be symmetric, complete, and satisfying triangle inequality. Further, let $c_{ij}$ (and $c (i, j)$) denote the Euclidean distance between vertices $i$ and $j$ in the graph. Without loss of generality, the vehicles are considered to be traveling at a unit speed. Hence, $c_{ij}$ denotes the cost of edge $(i, j) \in E$.

\subsection{Initial feasible solution}

The initial feasible solution is generated using a load balancing technique studied in \cite{min_max_vrp_LB_based_load_balancing} and subsequently utilized in the MD algorithm. In this technique, an assignment problem is formulated. Consider the depots of the vehicles, which are indexed by $j$, and the targets, which are indexed by $i$. Consider a binary variable $x_{ij}$ denoting whether target $i$ is allocated to depot $j$ or not. If target $i$ is allocated to depot $j$, then $x_{ij} = 1$, and if not, $x_{ij} = 0$. Let $c_{ij}$ denote the Euclidean distance between target $i$ and depot $j$. Consider the following integer program formulation:
\begin{align}
    & \min \sum_{i = 1}^{|T|} \sum_{j = 1}^{m} c_{ij} x_{ij} \label{eq: objective_function} \\
    \text{s.t.} \quad & \sum_{j = 1}^m x_{ij} = 1, \quad \forall i \in T, \label{eq: target_allocation} \\
    & \sum_{i = 1}^{|T|} x_{ij} = \begin{cases}
        \big{\lceil} \frac{|T|}{m} \big{\rceil}, j = 1, 2, \cdots, (|T| \ \mathrm{mod} \ m), \\
        \big{\lfloor} \frac{|T|}{m} \big{\rfloor}, j = (|T| \ \mathrm{mod} \ m) + 1, \cdots, m, \\
    \end{cases} \label{eq: vehicle_num_targets_assigned} \\
    & x_{ij} \in \{0, 1 \} \quad \forall i \in T \quad \forall j \in \{1, 2, \cdots, m\}. 
\end{align}
In the above formulation, constraint~\eqref{eq: target_allocation} ensures that each target is allocated to exactly one vehicle. Constraint~\eqref{eq: vehicle_num_targets_assigned} ensures that each vehicle is allocated approximately the same number of targets. Since the number of targets can be expressed as $|T| = p m + q,$ where $p$ and $q$ are integers, the first $q$ number of vehicles are allocated $p + 1 = \big{\lceil} \frac{|T|}{m} \big{\rceil}$ number of targets, and the other $m - q$ number of vehicles are allocated $p = \big{\lfloor} \frac{|T|}{m} \big{\rfloor}$ number of targets.
In the above formulation, the objective function in \eqref{eq: objective_function} allocates targets to depots that are in its vicinity. It should be noted that the above formulation can be solved as a linear program to obtain the optimal solution for the integer program.

It should be noted that if multiple vehicles start from the same depot location, then the depot location is perturbed about its initial location for each vehicle, similar to the MD algorithm \cite{MD_algorithm}. If the perturbation is not performed, varying target allocations between the vehicles that start at the same location results in no change in the objective function. Similar to \cite{MD_algorithm}, the initial locations were placed symmetrically on a circle of radius $0.1$ centered at the initial depot location (before perturbation). The angle of one of the perturbations was chosen randomly, and the perturbation angles for the other vehicles starting at the same depot were obtained such that perturbed depots are placed symmetrically on the circle.

\subsection{Local search}

\begin{figure*}[htb!]
    \centering
    \subfigure[Vehicle tours before removing target $t$]{\includegraphics[width = 0.3\textwidth]{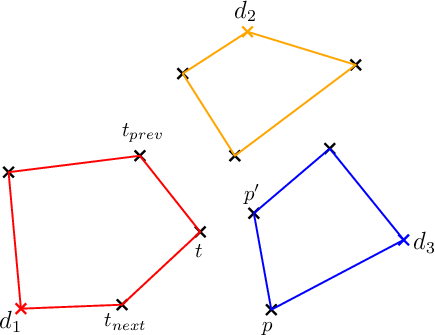}
    \label{subfig: vehicle_tours_before_removing_target_t}}
    \hfill
    \subfigure[Feasible tour for ``red" vehicle after removing target $t$]{\includegraphics[width = 0.3\textwidth]{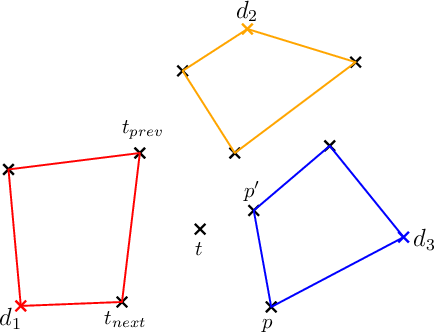}
    \label{subfig: vehicle_tour_after_removing_target_t}}
    \hfill
    \subfigure[Vehicle tour for ``blue" vehicle after inserting target $t$ between vertices $p$ and $p'$]{\includegraphics[width = 0.3\textwidth]{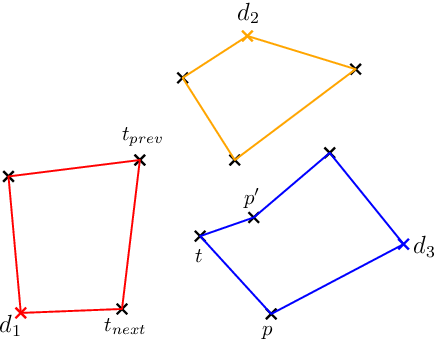}
    \label{subfig: vehicle_tour_after_inserting_target_t}}
    \caption{Vehicle tours after removing target $t$ from red vehicle and inserting into ``blue" vehicle}
    \label{fig: vehicle_tours_removing_t_inserting}
\end{figure*}

A local search is performed to improve the solution obtained. 
For this purpose, two neighborhoods are considered: Neighborhood $1$ and Neighborhood $2$. In both neighborhoods, a vehicle is first chosen.
\begin{itemize}
    \item In the first neighborhood (N1), a target is attempted to be removed from this vehicle and inserted into another vehicle. A depiction of N1 is shown in Fig.~\ref{fig: vehicle_tours_removing_t_inserting}.
    \item In the second neighborhood (N2), a target from the chosen vehicle is attempted to be swapped with a target allocated to another vehicle.
\end{itemize}
For the implementation of these neighborhoods, and to ensure that the computations are fast, three questions need to be answered:
\begin{itemize}
    \item How do we pick the vehicle to remove a target from (referred to as the maximal vehicle)?
    \item How do we sort the list of targets in the vehicle from which the target is attempted to be removed?
    \item How do we pick the vehicle in which a removed target needs to be inserted in (in N1) or swapped with (in N2)? (For reference, the vehicle picked for insertion in Fig.~\ref{fig: vehicle_tours_removing_t_inserting} is the ``blue" vehicle.)
\end{itemize}
To answer these questions, two proxy costs will be defined: a cost associated with removing a target from a vehicle, and a cost associated with inserting a target into a vehicle.

\subsubsection{Proxy cost for target removal}

Consider removing a target $t$ from the vehicle shown in red in Fig.~\ref{fig: vehicle_tours_removing_t_inserting}. It is desired to estimate the increase in the objective value of the red vehicle after removing target $t$ from it. Let $t_{prev}$ and $t_{next}$ denote the targets covered before and after target $t$ in the maximal vehicle's tour, and the maximal vehicle's tour be given by $(u_0 = D_j, u_1^j, u_2^j, \cdots, t_{prev}, t, t_{next}, \cdots, u_{n_j - 1}^j, u_{n_j}^j = D_j)$. It should be noted here that $D_j$ denotes the depot of the maximal vehicle. Hence, the estimated increase in cost associated with $t$ is defined to be
\begin{align} \label{eq: increase_in_objective_targ_removal}
\begin{split}
    &\text{increase in objective}_t \\
    &= \text{estimated new objective}_j - \text{old objective}_j \\
    &= J_j^{est} \left(d_{u_1^j}, d_{u_2^j}, \cdots, d_{t_{prev}}, d_{t_{next}}, \cdots, d_{u_{n_j - 1}^j} \right) \\
    & \quad\, - J_j \left(d_{u_1^j}, d_{u_2^j}, \cdots, d_{t_{prev}}, d_{t}, d_{t_{next}}, \cdots, d_{u_{n_j - 1}^j} \right),
\end{split}
\end{align}
where $d_v$ denotes the dwell time for a vertex $v$. Here, $J_j^{est} \left(d_{u_1^j}, d_{u_2^j}, \cdots, d_{t_{prev}}, d_{t_{next}}, \cdots, d_{u_{n_j - 1}^j} \right)$ denotes the estimated new objective for the $j\textsuperscript{th}$ vehicle after removing target $t$. For computing the estimated objective value, estimates for the dwell times of the other targets after removing target $t$, and the cost associated with the tour (TSP) after removing $t$, need to be obtained. To this end,
\begin{itemize}
    \item The dwell times of targets $u_1^j, u_2^j, \cdots, t_{prev}, t_{next}, \cdots, u_{n_j - 1}^j$ are maintained to be the same as before removing target $t$.
    \item Further, for the estimated tour cost, a ``savings" metric similar to the MD algorithm \cite{MD_algorithm} was considered. 
    The savings corresponding to target $t$ is given by
    \begin{align*}
        \text{savings}_t = c (t_{prev}, t) + c(t, t_{next}) - c(t_{prev}, t_{next}).
    \end{align*}
    Hence, the new tour cost is given by the previous tour cost minus the savings corresponding to target $t$.
\end{itemize}

\subsubsection{Proxy cost for target insertion}

Similar to the ``increase in objective" metric associated with removing a target, it is desired to define an ``increase in objective insert" metric associated with inserting target $t$ in the ``blue" vehicle shown in Fig.~\ref{fig: vehicle_tours_removing_t_inserting}. 
This metric estimates an increase in the objective value obtained with inserting target $t$ in the ``blue" vehicle (denoted with $i$ henceforth). For vehicle $i$,
\begin{itemize}
    \item The dwell times of the other targets covered by the vehicle, which are $u_1^i, \cdots, p, p', \cdots, u_{n_i - 1}^i$ are maintained to be the same as before inserting target $t$. Further, the dwell time of target $t$ before removing from vehicle $j$ and after inserting in vehicle $i$ are kept to be the same.
    \item Suppose $t$ is desired to be inserted between targets $p$ and $p'$ in the $i\textsuperscript{th}$ vehicle's tour. The insertion cost associated with inserting $t$ in vehicle $i$ between targets $p$ and $p'$ is defined to be
    \begin{align*}
        \text{cost increase}_{t, i, (p, p')} = c(p, t) + c(t, p') - c(p, p').
    \end{align*}
    This metric estimates the increase in the tour cost (TSP) with the insertion of $t$ between targets $p$ and $p'$. Since target $t$ can be inserted between any pair of vertices in vehicle $i$'s tour, the ``cost increase" for vehicle $i$ for inserting target $t$ is defined to be the minimum among all such costs. This choice of insertion will also yield the highest estimated increase in the objective value (as can be observed from \eqref{eq: objective_function_one_vehicle_given_partition}), since the discounting due to the exponential will be the least.
\end{itemize}
The ``increase in objective insert" corresponding to inserting target $t$ for vehicle $i$ can therefore be expressed as
\begin{align} \label{eq: increase_in_objective_insertion}
\begin{split}
    &\text{increase in objective insert}_{t, i} \\
    &= \text{estimated new objective}_i - \text{old objective}_i \\
    &= J_i^{est} \left(d_{u_1^i}, d_{u_2^i}, \cdots, d_t^{est}, \cdots, d_{u_{n_i - 1}^i} \right) \\
    & \quad\, - J_i \left(d_{u_1}^i, d_{u_2^i}, \cdots, d_{u_{n_i - 1}^i} \right).
\end{split}
\end{align}

\subsubsection{Analysis of proxy costs}

\begin{figure}[htb!]
    \centering
    \subfigure[Difference between real and proxy objective for removal of targets across all 43 instances]{
\includegraphics[width=.7\linewidth]{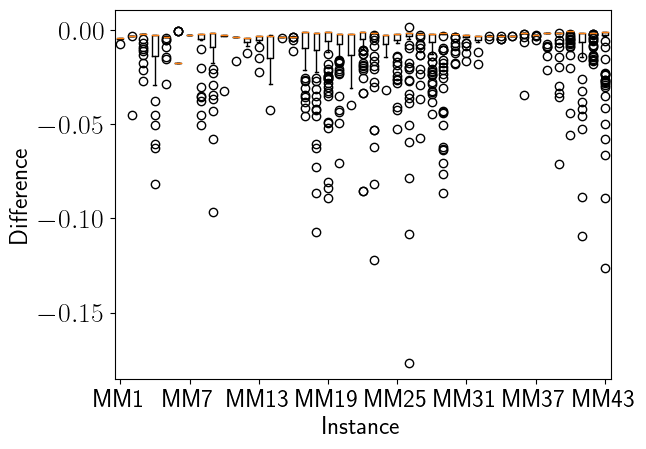}} 
    \subfigure[Comparison of real objective change and proxy objective change on MM30]{\includegraphics[width=.7\linewidth]{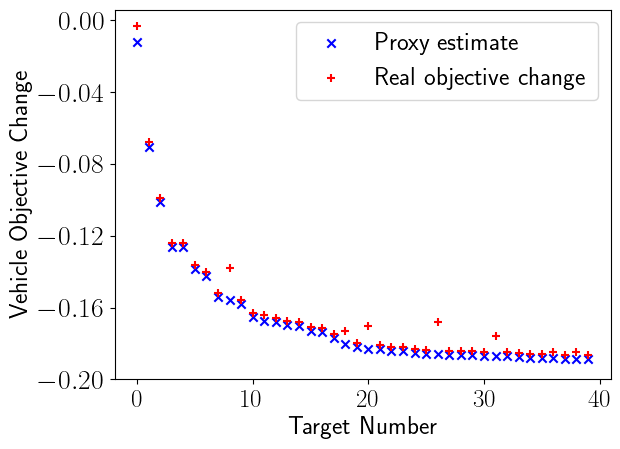}} 
    \caption{Proxy objective estimates compared to the real tour objective obtained through LKH and gradient descent}
    \label{fig:proxy_costs}
\end{figure}

To justify the use of our proxy objective estimates in local search, we compare the estimated change in vehicle objective with the real change. For each instance, the initial feasible solution based on the assignment problem is first generated. For each vehicle, we consider the removal of all targets. 
The results in Fig.~\ref{fig:proxy_costs} show the difference between the true objective value and the proxy estimates of each vehicle over 43 instances taken from \cite{MD_algorithm}. It should be noted that the instances considered are diverse with respect to the number of targets, the number of vehicles, and the probability distribution from which the target locations were obtained. From this figure, we can observe that the proxy costs provide a relatively accurate estimate of the objective value with target removal. We note that though all values in the box plot should be non-positive, we see a few targets where the proxy estimate is larger than the real objective, which arises due to rounding of edge costs for running LKH.

\subsubsection{Variations of vehicle choice considered for removal of target}

\begin{figure}[htb!]
    \centering
    \subfigure[$\tau = 2$]{\includegraphics[width = 0.3\textwidth]{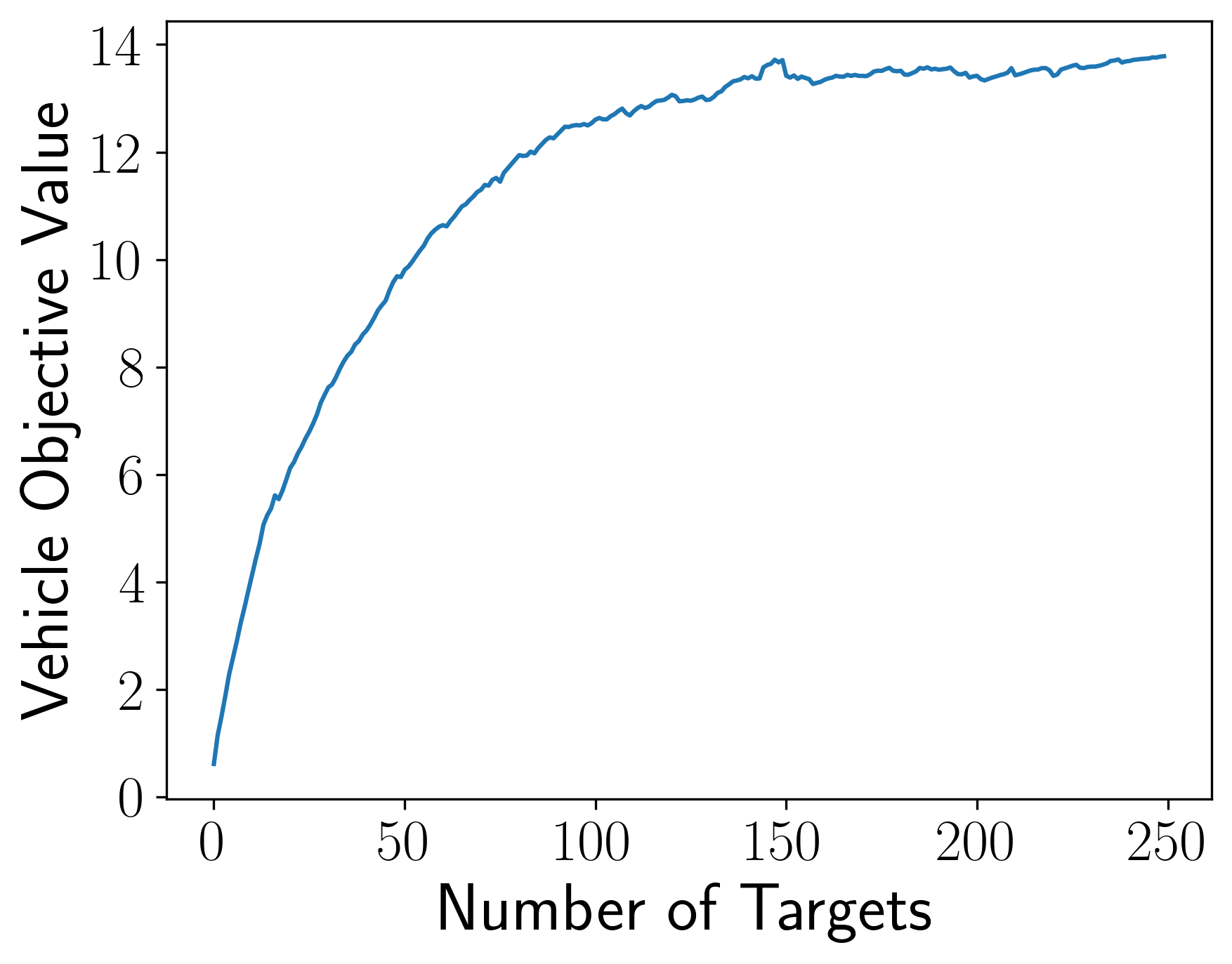}}
    \subfigure[$\tau = 20$]{\includegraphics[width = 0.3\textwidth]{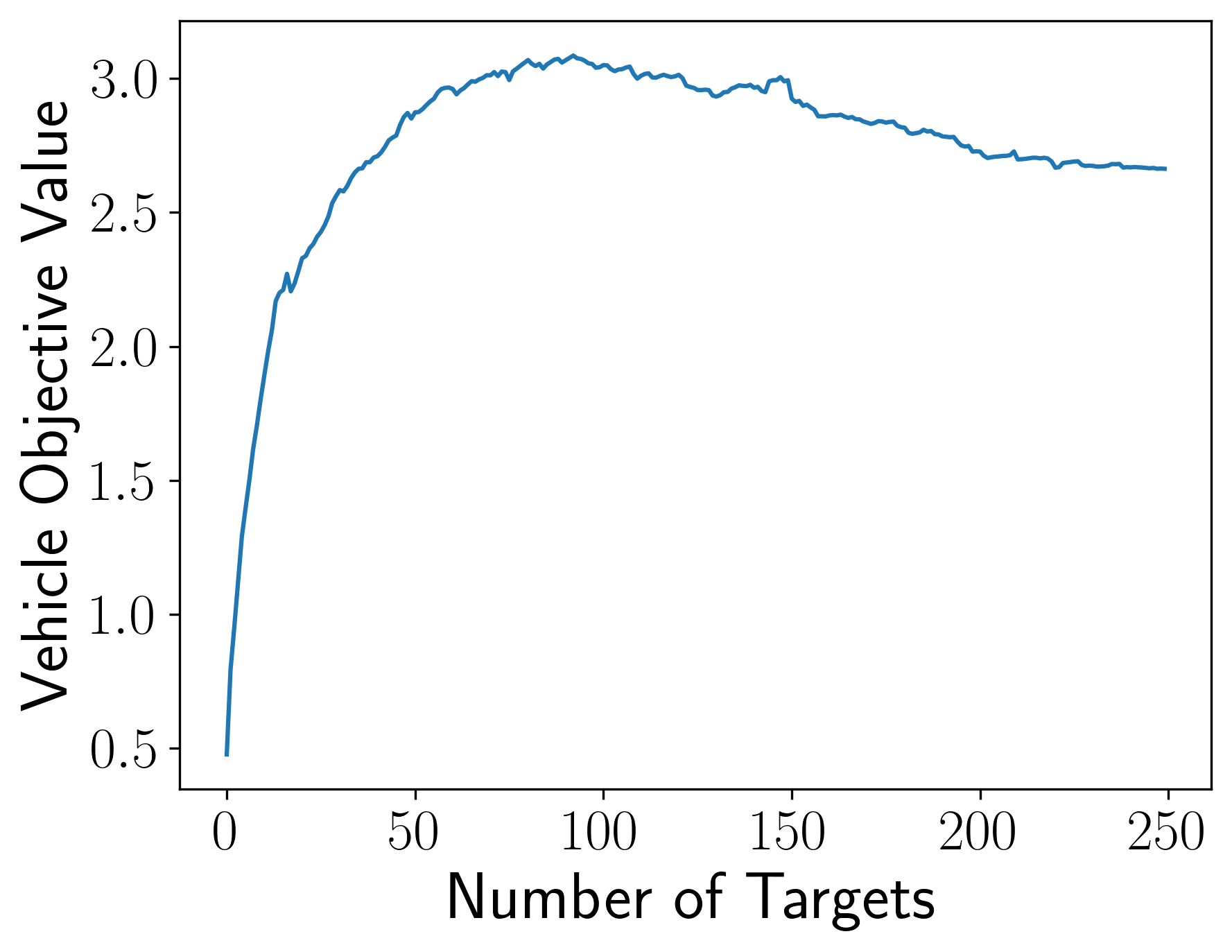}}
    \caption{Sample vehicle objective value variation with number of targets for instance MM22 (from \cite{MD_algorithm}) for $\alpha = 8.77\cdot 10^{-4}, \tau = 2$ and $\alpha = 8.77\cdot 10^{-4}, \tau = 20$}
    \label{fig:targets}
\end{figure}

Having defined the estimated increase in the objective value associated with removing and inserting a target, it is desired to identify the vehicle (termed the maximal vehicle) from which a target is attempted to be removed (the ``red" vehicle in Fig.~\ref{fig: vehicle_tours_removing_t_inserting}). While the choice of the maximal vehicle is immediate in a min-max problem, which is the vehicle that has the highest tour cost, such a choice of maximal vehicle is not immediate for our considered problem. To this end, four variations in the decision of the vehicle to be picked were considered. While one of the variations described below will be a natural choice for this problem, the justification for three of the variations is drawn using Fig.~\ref{fig:targets}. In this figure, the variation in the objective value of the first vehicle, obtained from \eqref{eq: objective_function_one_vehicle_given_partition}, with changing number of targets is shown for instance MM22 (taken from \cite{MD_algorithm}).
The value for $\alpha$ was computed as $\frac{1}{TSP},$ where $TSP$ denotes the tour cost to cover all targets by a single vehicle, which was obtained using LKH.\footnotemark Further, two variations in $\tau$ were considered for this plot.
\footnotetext{The reason for considering such an $\alpha$ value based on the TSP cost was to ensure that the discounting factor corresponding to the exponential in \eqref{eq: objective_function_one_vehicle_given_partition} is not very large. In such cases, it was observed that one of the vehicles covers all targets, whereas other vehicles cover no targets or one target.}

The variations considered in the choice of maximal vehicle and the reason behind considering each variation are as follows:
\begin{itemize}
    \item Vehicle with minimum objective value: The intuition behind such a vehicle choice is that a vehicle with a low objective value will be indicative of a vehicle with a high tour cost. This, in turn, would lead to a high discounting value due to the exponential in the objective function (refer to equation~\eqref{eq: objective_function_single_vehicle}). 
    \item Vehicle with maximum objective value: The intuition behind removing a target from such a vehicle is that the considered vehicle's objective value would not significantly be affected by removing a target, as shown in Fig.~\ref{fig:targets}. However, by inserting the removed target in another vehicle, a net improvement in the objective value can be made.
    \item Vehicle with the maximum number of targets: If a vehicle visits many targets, then the removal of a target can lead to an increase in the objective value due to an increase in the exponential in the objective function (refer to equation~\eqref{eq: objective_function_single_vehicle} and Fig.~\ref{fig:targets}).
    \item Vehicle with maximum tour cost: The reason for this choice is the same as the choice of vehicle with the maximum number of targets. However, this variation is considered since a vehicle with the maximum tour cost need not necessarily be the same as the vehicle with the maximum number of targets.
\end{itemize}

It is desired to determine which of the variations will yield a good improvement with respect to the initial solution. To this end, it is desired to study the influence of vehicle choice on the objective value obtained after a local search for different instances. For this purpose, we will first describe neighborhood $1$, which will later be used for this study (and for the heuristic).

\subsubsection{Neighborhood 1: 1 pt. move}

For a chosen vehicle (``maximal" vehicle) from which a target is desired to be removed, an order of considering targets in the maximal vehicle is first desired to be constructed. To this end, the ``increase in objective" metric associated with removing a target, which was previously defined, is used. Since it is desired to remove a target that can provide a high increase in the objective value, the targets in the maximal vehicle are sorted in the decreasing order of this metric. 

Suppose target $t$ is considered to be removed from the maximal vehicle. The vehicle in which the target $t$ will be inserted needs to be obtained. To this end, the ``increase in objective insert" metric previously defined is used. Among all vehicles other than the maximal vehicle, the vehicle with the highest increase in objective metric is chosen for inserting target $t$ (since a maximization problem is considered). Hence, if the estimated objective value corresponding to the new solution obtained after removing $t$ from the maximal vehicle and inserting it in another vehicle is higher than the previous objective value, then the obtained solution is chosen as the new incumbent solution. Further, LKH is used to construct the tours of the maximal vehicle and the vehicle considered for insertion, and gradient descent is run to obtain the dwell time for each target covered by the two vehicles. If a better solution was not obtained by removing target $t$, the next target in the sorted list by the ``increase in objective" metric defined is considered from the maximal vehicle. The same steps are performed till a better solution is obtained or all targets from the maximal vehicle are considered.

\subsubsection{Impact of maximal vehicle considered for the local search}

\begin{figure}[htb!]
    \centering
    \includegraphics[width = .9\linewidth]{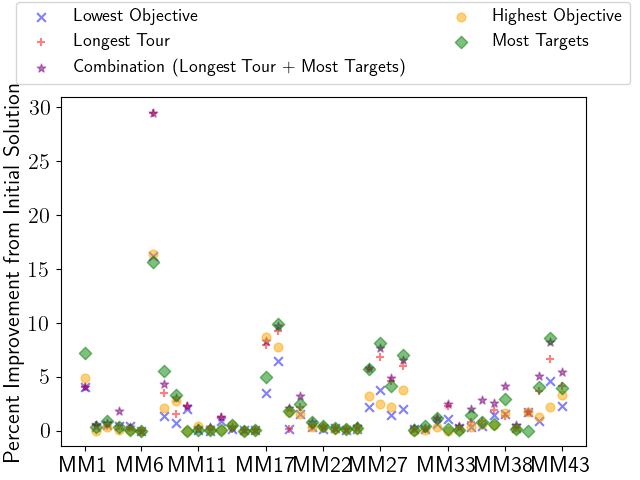}
    \caption{Percentage improvement using 1 pt. move by heuristic for ``vehicle to remove target from" based on different vehicle choices}
    \label{fig:percentage_improvement}
\end{figure}

For the neighborhood previously considered, a maximal vehicle needs to be identified, which is the vehicle from which a target is attempted to be removed. For this purpose, a total of four variations were considered, as previously mentioned. The percentage improvement from the local search using neighborhood $1$ with respect to the initial solution for all four methods are reported for a total of $43$ instances, which were taken from \cite{MD_algorithm}. The obtained percentage improvement for the four variations are shown in Fig.~\ref{fig:percentage_improvement}, and the minimum, median, mean, and maximum percentage improvements are summarized in Table~\ref{tab:percentage_improvement}. It can be observed that the variations corresponding to choosing the vehicle with the most number of targets and the vehicle with the highest tour cost yielded the best improvement with respect to the initial solution.

Since both variations yield a good improvement with respect to the initial solution and different improvement values, the previous neighborhood considered was modified to consider two types of ``maximal" vehicle. First, the maximal vehicle is considered to be the vehicle with the longest tour. If a better solution was obtained, then the solution obtained was considered to be the new incumbent solution. If no improvement could be obtained, then the same neighborhood search is performed, but with the maximal vehicle considered to be the vehicle with the most number of targets. The percentage improvement with this modified neighborhood is also reported in Fig.~\ref{fig:percentage_improvement} and Table~\ref{tab:percentage_improvement}, and can be observed to provide a better solution than the solutions previously obtained. Hence, the modified neighborhood considering both types of maximal vehicles would be considered for the heuristic.

\begin{table}[htb!]
	\centering
         \caption{Percentage improvement of Local Search by heuristic for \lq vehicle to remove target from\rq}
        \label{tab:percentage_improvement}
	\begin{tabular}{|c|c|c|c|c|}
	\hline
	\textbf{Method} & \textbf{Min.} & \textbf{Median} & \textbf{Mean} & \textbf{Max.} \\\hline
     \textit\textbf{{Lowest Objective}} & 0 & 0.5 & 1.5 & 16.2\\\hline
     \textit\textbf{{Highest Objective}} & 0 & 0.33 & 1.68 & 16.37\\\hline
     \textit\textbf{{Most Targets}} & 0 & 0.61 & 2.43 & 15.64\\\hline
	\textit\textbf{{Longest Tour}} & 0 & 0.88 & 2.64 & 29.47\\\hline
	Combination (Longest Tour & \multirow{2}{*}{0} & \multirow{2}{*}{1.86} & \multirow{2}{*}{3.16} & \multirow{2}{*}{29.47} \\
          + Most Targets) & & & & \\ \hline
	\end{tabular}
\end{table}

\subsubsection{Neighborhood 2: 1 pt. swap}

Using the identified combination of maximal vehicles using neighborhood 1, the same combination of maximal vehicles is used for neighborhood 2. In this neighborhood, it is desired to select a target from a maximal vehicle and swap it with a target from another vehicle to obtain a better solution. Using the observation of the combination of maximal vehicles previously considered, first
\begin{itemize}
    \item The vehicle corresponding to the longest tour is utilized as the maximal vehicle, and neighborhood 2 is used.
    \item If a better solution was not obtained, the vehicle corresponding to the maximum number of targets is utilized as the maximal vehicle, and neighborhood 2 is used.
\end{itemize}
Having chosen the maximal vehicle (either the vehicle with the longest tour or the maximum number of targets), the targets are ordered in the maximal vehicle in the descending order of the increase in objective metric given in \eqref{eq: increase_in_objective_targ_removal}. Suppose target $t$ is removed from the maximal vehicle. Similar to neighborhood 1, the target is inserted in the vehicle with the highest increase in the objective value corresponding to target $t$, which is defined in \eqref{eq: increase_in_objective_insertion}. Let the vehicle in which $t$ is inserted be denoted by $i$. It is now desired to remove a target from vehicle $i$ and insert it into the maximal vehicle.

Similar to target removal from the maximal vehicle, the targets in vehicle $i$ are ordered in the decreasing order of the increase in objective value given in \eqref{eq: increase_in_objective_targ_removal}. However, it should be noted that contrary to target removal from the maximal vehicle, the objective value of vehicle $i$ before removing a target is an estimate. That is, in \eqref{eq: increase_in_objective_targ_removal}, ``old objective$_i$" is the estimated objective value of vehicle $i$ obtained after inserting $t$. Having ordered the targets in vehicle $i$ in the decreasing order of the estimated increase in objective, each of the targets is attempted to be inserted in the maximal vehicle. Suppose $t_i$ is attempted to be inserted in the maximal vehicle's tour. For this purpose, the increase in the objective value corresponding to the insertion of $t_i$ is obtained for the maximal vehicle, using \eqref{eq: increase_in_objective_insertion}. 
If a better solution is obtained, LKH is used to optimize the vehicle tours, and the dwell times are obtained using gradient descent for the two vehicles. The same steps are performed till a better solution is obtained or all targets from the maximal vehicle are considered.

\subsection{Perturbation of solution}

\begin{figure}[htb!]
    \centering
    \includegraphics[width = 0.45\linewidth]{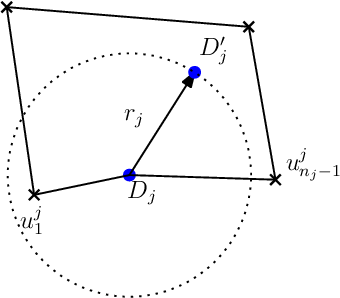}
    \caption{Depiction of perturbation for the $j\textsuperscript{th}$ vehicle}
    \label{fig: perturbation_j_vehicle}
\end{figure}

The solution obtained from local search is a local minimum, and cannot be improved using the neighborhoods defined. Hence, it is desired to break from the local minimum. To this end, the solution is perturbed similar to the MD algorithm \cite{MD_algorithm} by perturbing the depot locations. For this purpose, for the $j\textsuperscript{th}$ vehicle, the average distance $r_j$ of the depot $D_j$ from the two targets connected to it in the vehicle's tour in the current solution is obtained for $j = 1, 2, \cdots, m.$ Then, in the first round of perturbation, each depot location is perturbed from its initial location by a random angle to a location that is at a distance $r_j$ from its initial location, as shown in Fig.~\ref{fig: perturbation_j_vehicle}. Using the same allocation of targets and the corresponding dwell times from the solution obtained from the local search for each vehicle, a feasible solution is generated for each vehicle. A local search is then performed on the graph with the perturbed depot locations to obtain a new allocation of targets for each vehicle and corresponding dwell times. Using this new solution in the original graph, wherein the depots are restored to their initial locations, a local search is then performed to obtain a new potential solution. If the obtained solution is better than the incumbent solution, then the incumbent solution is updated and the perturbation step is restarted. 

Similar to \cite{MD_algorithm}, the perturbation step is performed for five consecutive iterations till no improvement is obtained. The perturbation angle for each depot is set to be $144^\circ$ from the previous iteration's perturbation angle. Hence, the sixth perturbation will be the same as the first perturbation.  

\section{Simulation Results}

The results we present in this section were generated using Python 3.8 and Julia 1.9.2 on a laptop with AMD Ryzen 5 4600H running at 3 GHz with 8 GB RAM. We first present results corresponding to the single vehicle case, which can be solved to optimality. The multi-vehicle case results are then presented, wherein the performance of different variations of the presented heuristic are analyzed on a diverse set of instances.

\subsection{Single vehicle results}

The single-vehicle case was solved to optimality in two steps: solving the integer program corresponding to the single-vehicle TSP using branch and cut, and obtaining the dwell times using gradient descent. To this end, branch and cut was implemented in Julia using the JuMP package and was solved using Gurobi \cite{gurobi}. The integer program was warm-started using the solution from LKH \cite{LKH_TSP_solver}, a state-of-the-art heuristic for single-vehicle TSP. Further, for the branch and cut implementation, subtour elimination constraints were used as cutting planes and lazy cuts using lazy callback functions.

The instances considered were taken from standard US datasets from TSPLIB \cite{TSPLIB95}. For each of the instances, three values of $\alpha$ were considered. The $\alpha$ values chosen were equal to $\frac{0.5}{TSP}, \frac{1}{TSP}, \frac{2}{TSP},$ where $TSP$ here denotes the tour cost associated with covering all the targets in the graph, which was obtained using LKH. The reason for the choice of such $\alpha$'s were to ensure that the discounting associated with the exponential terms in the objective function are not very high.

\textbf{Remark:} When the $\alpha$ was set to larger values, i.e., when the discounting due to the exponential is high, it was observed that the optimal solution for three vehicle problems with ten targets were such that one of the vehicles covers a majority of the targets, and the other vehicles covers at most one target.


For the instances, three values of $\tau$ were considered: $0.5, 1,$ and $2.$ The results obtained for the single vehicle case are summarized in Table~\ref{tab: single_vehicle_results}. From this table, we can observe that
\begin{enumerate}
    \item The optimal solution could be obtained within about 30 seconds for all instances.
    \item The dwell time for each target is approximately the same due to the same $\tau$ value for all targets.
    \item Increasing $\tau$ increases the dwell time for each target, whereas increasing $\alpha$ decreases the dwell time.
\end{enumerate}

\begin{table*}[htb!]
\centering
\caption{Results for single vehicle case}
\label{tab: single_vehicle_results}	\begin{tabular}{|c|c|c|c|c|c|c|c|c|c|c|}
	\hline
\multirow{2}{*}{\textbf{Instance}} &\textbf{No. of} &\multirow{2}{*}{\textbf{TSP cost}$^*$} &\multirow{2}{*}{$\alpha$ (s\textsuperscript{-1})} &\multirow{2}{*}{$\tau$ (s)} &\multicolumn{4}{c|}{\textbf{Dwell Time} (s)}&\textbf{Branch and cut} &\textbf{Gradient}\\
\cline{6-9}
& {\textbf{targets}} & & & & \textbf{Min}& \textbf{Median}& \textbf{Mean}& \textbf{Max} & \textbf{time} (s) & \textbf{descent time} (s) \\

\hline
 	\multirow{9}{*}{rd100} & \multirow{9}{*}{99} & \multirow{9}{*}{7910.4} & \multirow{3}{*}{6.37e-05} & 0.5 & 11.37 & 11.37 & 11.37 & 11.38 & 1.09 & 1.50e-02 \\
	& & & & 1.0 & 16.97 & 17.03 & 17.01 & 17.03 & 1.19 & 9.00e-03 \\
	& & & & 2.0 & 24.4 & 24.49 & 24.55 & 24.64 & 1.1 & 1.10e-02 \\

\cline{4-11}
	& & &\multirow{3}{*}{1.27e-04} & 0.5 & 8.76 & 8.77 & 8.78 & 8.78 & 1.08 & 5.00e-03 \\
	& & & & 1.0 & 12.65 & 12.75 & 12.72 & 12.76 & 1.11 & 6.00e-03 \\
	& & & & 2.0 & 17.56 & 17.82 & 17.71 & 17.85 & 1.13 & 8.00e-03 \\

\cline{4-11}
	& & &\multirow{3}{*}{2.55e-04} & 0.5 & 6.45 & 6.47 & 6.48 & 6.5 & 1.13 & 3.00e-03 \\
	& & & & 1.0 & 8.96 & 9.09 & 9.04 & 9.11 & 1.12 & 3.00e-03 \\
	& & & & 2.0 & 12.01 & 12.06 & 12.1 & 12.22 & 1.06 & 4.00e-03 \\

\hline
 	\multirow{9}{*}{bier127} & \multirow{9}{*}{126} & \multirow{9}{*}{118293.52} & \multirow{3}{*}{4.27e-06} & 0.5 & 19.16 & 19.16 & 19.16 & 19.17 & 1.04 & 2.10e-02 \\
	& & & & 1.0 & 30.45 & 30.45 & 30.45 & 30.46 & 1.05 & 3.10e-02 \\
	& & & & 2.0 & 47.1 & 47.16 & 47.17 & 47.19 & 1.04 & 4.50e-02 \\

\cline{4-11}
	& & &\multirow{3}{*}{8.54e-06} & 0.5 & 17.52 & 17.52 & 17.52 & 17.52 & 1.06 & 1.80e-02 \\
	& & & & 1.0 & 27.57 & 27.58 & 27.58 & 27.58 & 1.04 & 2.30e-02 \\
	& & & & 2.0 & 42.14 & 42.21 & 42.22 & 42.25 & 1.05 & 4.00e-02 \\

\cline{4-11}
	& & &\multirow{3}{*}{1.71e-05} & 0.5 & 15.3 & 15.3 & 15.31 & 15.31 & 1.05 & 1.40e-02 \\
	& & & & 1.0 & 23.7 & 23.72 & 23.72 & 23.72 & 1.06 & 1.90e-02 \\
	& & & & 2.0 & 35.57 & 35.64 & 35.67 & 35.72 & 1.05 & 2.50e-02 \\

\hline
 	\multirow{9}{*}{pr152} & \multirow{9}{*}{151} & \multirow{9}{*}{73683.64} & \multirow{3}{*}{6.79e-06} & 0.5 & 17.37 & 17.37 & 17.37 & 17.38 & 49.27 & 2.20e-02 \\
	& & & & 1.0 & 27.31 & 27.31 & 27.32 & 27.32 & 51.26 & 3.90e-02 \\
	& & & & 2.0 & 41.68 & 41.8 & 41.77 & 41.81 & 49.28 & 5.00e-02 \\

\cline{4-11}
	& & &\multirow{3}{*}{1.36e-05} & 0.5 & 15.3 & 15.31 & 15.3 & 15.31 & 49.38 & 1.80e-02 \\
	& & & & 1.0 & 23.69 & 23.72 & 23.71 & 23.72 & 49.57 & 4.40e-02 \\
	& & & & 2.0 & 35.54 & 35.64 & 35.65 & 35.69 & 52.36 & 4.10e-02 \\

\cline{4-11}
	& & &\multirow{3}{*}{2.72e-05} & 0.5 & 12.83 & 12.84 & 12.84 & 12.84 & 53.94 & 1.50e-02 \\
	& & & & 1.0 & 19.45 & 19.49 & 19.48 & 19.5 & 53.34 & 3.10e-02 \\
	& & & & 2.0 & 28.45 & 28.5 & 28.59 & 28.66 & 50.7 & 2.60e-02 \\

\hline
 	\multirow{9}{*}{d198} & \multirow{9}{*}{197} & \multirow{9}{*}{15808.65} & \multirow{3}{*}{3.17e-05} & 0.5 & 11.15 & 11.16 & 11.16 & 11.16 & 11.07 & 1.70e-02 \\
	& & & & 1.0 & 16.61 & 16.68 & 16.66 & 16.68 & 11.24 & 3.50e-02 \\
	& & & & 2.0 & 23.83 & 23.84 & 23.99 & 24.08 & 11.56 & 2.70e-02 \\

\cline{4-11}
	& & &\multirow{3}{*}{6.34e-05} & 0.5 & 8.67 & 8.7 & 8.69 & 8.7 & 11.07 & 1.20e-02 \\
	& & & & 1.0 & 12.5 & 12.56 & 12.58 & 12.62 & 11.31 & 1.60e-02 \\
	& & & & 2.0 & 17.34 & 17.36 & 17.5 & 17.62 & 11.02 & 1.90e-02 \\

\cline{4-11}
	& & &\multirow{3}{*}{1.27e-04} & 0.5 & 6.42 & 6.45 & 6.45 & 6.47 & 11.48 & 9.00e-03 \\
	& & & & 1.0 & 8.93 & 8.99 & 9.01 & 9.07 & 11.53 & 1.10e-02 \\
	& & & & 2.0 & 11.96 & 12.16 & 12.03 & 12.17 & 11.38 & 1.40e-02 \\

\hline
 	\multirow{9}{*}{pr226} & \multirow{9}{*}{225} & \multirow{9}{*}{80370.26} & \multirow{3}{*}{6.23e-06} & 0.5 & 15.95 & 15.95 & 15.95 & 15.95 & 34.36 & 4.40e-02 \\
	& & & & 1.0 & 24.82 & 24.82 & 24.84 & 24.84 & 32.13 & 6.10e-02 \\
	& & & & 2.0 & 37.45 & 37.49 & 37.56 & 37.6 & 32.17 & 7.90e-02 \\

\cline{4-11}
	& & &\multirow{3}{*}{1.25e-05} & 0.5 & 13.89 & 13.89 & 13.89 & 13.9 & 32.8 & 4.00e-02 \\
	& & & & 1.0 & 21.25 & 21.3 & 21.29 & 21.3 & 30.36 & 3.70e-02 \\
	& & & & 2.0 & 31.46 & 31.47 & 31.59 & 31.66 & 35.14 & 5.30e-02 \\

\cline{4-11}
	& & &\multirow{3}{*}{2.49e-05} & 0.5 & 11.49 & 11.5 & 11.49 & 11.5 & 34.81 & 2.20e-02 \\
	& & & & 1.0 & 17.16 & 17.23 & 17.21 & 17.23 & 33.04 & 4.50e-02 \\
	& & & & 2.0 & 24.73 & 24.78 & 24.88 & 24.98 & 33.16 & 5.80e-02 \\

\hline
 	\multirow{9}{*}{gr229} & \multirow{9}{*}{228} & \multirow{9}{*}{1635.31$^{\dagger}$} & \multirow{3}{*}{3.06e-04} & 0.5 & 3.77 & 3.79 & 3.8 & 3.84 & 7.2 & 7.00e-03 \\
	& & & & 1.0 & 4.96 & 4.98 & 4.98 & 5.04 & 7.13 & 8.00e-03 \\
	& & & & 2.0 & 6.24 & 6.36 & 6.37 & 6.75 & 7.13 & 9.00e-03 \\

\cline{4-11}
	& & &\multirow{3}{*}{6.13e-04} & 0.5 & 2.5 & 2.5 & 2.51 & 2.54 & 7.25 & 4.00e-03 \\
	& & & & 1.0 & 3.12 & 3.15 & 3.16 & 3.21 & 7.34 & 3.50e-02 \\
	& & & & 2.0 & 3.7 & 3.82 & 3.8 & 3.83 & 7.81 & 1.40e-02 \\

\cline{4-11}
	& & &\multirow{3}{*}{1.23e-03} & 0.5 & 1.58 & 1.59 & 1.59 & 1.6 & 7.02 & 9.00e-03 \\
	& & & & 1.0 & 1.88 & 1.89 & 1.91 & 1.92 & 6.92 & 7.00e-03 \\
	& & & & 2.0 & 2.21 & 2.25 & 2.23 & 2.26 & 7.19 & 2.40e-02 \\
	\hline
 \multicolumn{11}{l}{$^*$: If we calculate these by rounding each edge weight, we get results consistent with the known optimal solution for these instances.} \\
 \multicolumn{11}{l}{$^\dagger$: The known optimal solution corresponds to an edge cost computation different from Euclidean distance.}
	\end{tabular}
\end{table*}

\subsection{Multi-vehicle case}

The instances considered in this study were taken from the MD datasets \cite{MD_algorithm}. For each instance, the heuristic proposed in Section~\ref{sect: heuristic} was utilized to obtain a feasible solution. 
We note here that for the heuristic, we choose to use the combination of the vehicle with the longest tour cost and the vehicle with the most number of targets for the maximal vehicle, i.e., to remove a target from. This is based on the results of our first parametric study, which is summarized in Fig.~\ref{fig:percentage_improvement}. Additionally, we perform a parametric study on the heuristic, wherein we compare the effects of the neighborhoods considered and the number of maximal vehicles considered on our final objective values as well as the computation times. In particular, we compare the following variations:
\begin{itemize}
\item Neighborhood Type
\begin{itemize}
    \item One point move
    \item One point move and One point swap
\end{itemize}
\item Number of vehicles to consider from heuristic
\begin{itemize}
    \item One maximal vehicle based on the longest tour first, followed by a maximal vehicle based on the most number of targets.
    \item Top two vehicles with the two highest longest tours, followed by the two vehicles with the most number of targets covered.
\end{itemize}
\end{itemize}
In addition, to demonstrate the efficacy of the heuristic on different $\alpha$ and $\tau$ values, four variations were considered for each instance, which are
\begin{itemize}
    \item $\alpha = \frac{1}{TSP}, \tau = 1,$
    \item $\alpha = \frac{1}{TSP}, \tau = 2,$
    \item $\alpha = \frac{2}{TSP}, \tau = 1,$ and
    \item $\alpha = \frac{2}{TSP}, \tau = 1.$
\end{itemize}
Here, $TSP$ denotes the tour cost corresponding to a single vehicle covering all the targets, which is obtained using LKH.

\begin{figure}[htb!]
    \centering
    \subfigure[Percent improvement from initial solution for each local search mode by instance]{\includegraphics[width = 0.4\textwidth]{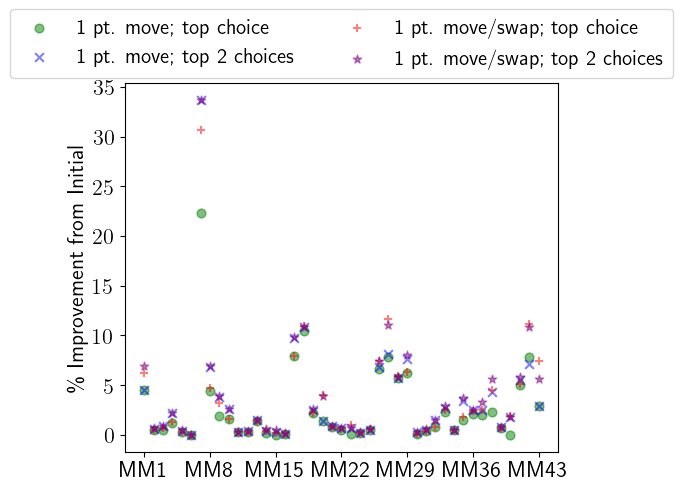}} \\
    \subfigure[Compute times for each local search mode by instance]{\includegraphics[width = 0.4\textwidth]{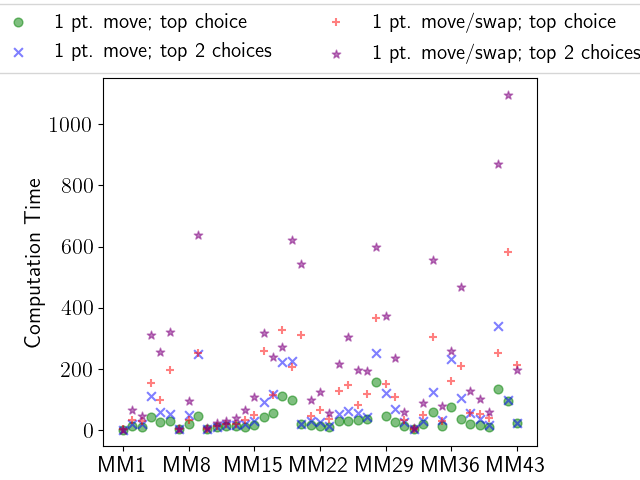}}
    \caption{Comparison between neighborhood search parameters with $\alpha=\frac{1}{TSP}$, $\tau=1$}
    \label{fig:local_search_parametric_comp11}
\end{figure}

\begin{table*}[htb!]
	\centering
    \caption{Percentage Improvement from Initial Solution and Computation Time by neighborhood search parameters}
    \label{tab:multivehicle_results}
	\begin{tabular}{|c|c|c|c|c|c|c|c|}
	\hline
	\textbf{$\alpha$} &\textbf{$\tau$}&\textit{\textbf{Value}} &\textbf{Method} & \textbf{Minimum} & \textbf{Median} & \textbf{Mean} & \textbf{Maximum} \\\hline
 \multirow{8}{*}{$\frac1{TSP}$}&\multirow{8}{*}{1}&&1 pt. move; top choice &  9.45e-06 & 1.43 & 3.05 & 29.7 \\\cline{4-8}
	&&\textit{\textbf{Percentage}}&\textit\textbf{{1 pt. move; top 2 choices}} &  9.45e-06 & 1.91 & 3.69 & 34.04 \\\cline{4-8}
    &&\textit{\textbf{Improvement}}&1 pt. move/swap; top choice &9.45e-06 & 1.56 & 3.42 & 30.68 \\\cline{4-8}
	&&&\textit\textbf{{1 pt. move/swap; top 2 choices}} & 9.45e-06 & 2.2 & 3.99 & 34.04\\\cline{3-8}\cline{3-8}\cline{3-8}\cline{3-8}
 
	&&&\textit\textbf{{1 pt. move; top choice}} & 1.78 & 31.85 & 45.88 & 171.93 \\\cline{4-8}
        &&\textit{\textbf{Computation}}&\textit\textbf{{1 pt. move; top 2 choices}} & 1.87 & 52.57 & 96.83 & 414.78 \\\cline{4-8} 
	&&\textit{\textbf{Time (s)}}&\textit\textbf{{1 pt. move/swap; top choice}} &  4.58 & 88.34 & 140.87 & 504.57 \\\cline{4-8}
	&&&\textit\textbf{{1 pt. move/swap; top 2 choices}} & 2.31 & 165.29 & 274.21 & 1106.77  \\\hline
\multirow{8}{*}{$\frac1{TSP}$}&\multirow{8}{*}{2}&&1 pt. move; top choice & 0.00e+00 & 1.21 & 2.98 & 29.76 \\\cline{4-8}
	&&\textit{\textbf{Percentage}}&\textit\textbf{{1 pt. move; top 2 choices}} & 9.45e-06 & 2.20 & 3.67 & 34.08\\\cline{4-8}
    &&\textit{\textbf{Improvement}}&1 pt. move/swap; top choice & 9.45e-06 & 1.41 & 3.36 & 30.73\\\cline{4-8}
	&&&\textit\textbf{{1 pt. move/swap; top 2 choices}} & 9.45e-06 & 2.20 &  3.91 & 34.08\\\cline{3-8}\cline{3-8}\cline{3-8}\cline{3-8}
 
	&&&\textit\textbf{{1 pt. move; top choice}} & 1.82 & 29.99 & 40.56 & 145.35 \\\cline{4-8}
        &&\textit{\textbf{Computation}}&\textit\textbf{{1 pt. move; top 2 choices}} & 1.88 & 52.47 & 94.74 & 428.38 \\\cline{4-8} 
	&&\textit{\textbf{Time (s)}}&\textit\textbf{{1 pt. move/swap; top choice}} & 3.94 & 86.93 & 144.46 & 604.44 \\\cline{4-8}
	&&&\textit\textbf{{1 pt. move/swap; top 2 choices}} & 2.30 & 140.59 & 281.96 & 1241.92 \\\hline
 \multirow{8}{*}{$\frac2{TSP}$}&\multirow{8}{*}{1}&&1 pt. move; top choice &  0.00e+00 & 2.83 & 6.38 & 68.73  \\\cline{4-8}
	&&\textit{\textbf{Percentage}}&\textit\textbf{{1 pt. move; top 2 choices}} & 1.89e-05 & 3.85 & 7.75 & 79.8 \\\cline{4-8}
    &&\textit{\textbf{Improvement}}&1 pt. move/swap; top choice & 1.89e-05 & 3.16 & 7.12 & 70.61 \\\cline{4-8}
	&&&\textit\textbf{{1 pt. move/swap; top 2 choices}} &  1.89e-05 & 4.43 & 8.48 & 79.8 \\\cline{3-8}\cline{3-8}\cline{3-8}\cline{3-8}
 
	&&&\textit\textbf{{1 pt. move; top choice}} &  2.12 & 30.54 & 43.04 & 178.28  \\\cline{4-8}
        &&\textit{\textbf{Computation}}&\textit\textbf{{1 pt. move; top 2 choices}} & 2.05 & 53.95 & 95.08 & 390.1 \\\cline{4-8} 
	&&\textit{\textbf{Time (s)}}&\textit\textbf{{1 pt. move/swap; top choice}} & 2.05 & 79.0 & 137.84 & 551.94 \\\cline{4-8}
	&&&\textit\textbf{{1 pt. move/swap; top 2 choices}} & 2.51 & 153.99 & 272.83 & 972.27 \\\hline
 \multirow{8}{*}{$\frac2{TSP}$}&\multirow{8}{*}{2}&&1 pt. move; top choice &   0.00e+00 & 2.81 & 6.24 & 68.94 \\\cline{4-8}
	&&\textit{\textbf{Percentage}}&\textit\textbf{{1 pt. move; top 2 choices}} & 1.89e-05 & 4.33 & 7.84 & 80.21\\\cline{4-8}
    &&\textit{\textbf{Improvement}}&1 pt. move/swap; top choice & 1.89e-05 & 3.03 & 7.27 & 70.75 \\\cline{4-8}
	&&&\textit\textbf{{1 pt. move/swap; top 2 choices}} &  1.89e-05 & 4.46 & 8.67 & 80.21 \\\cline{3-8}\cline{3-8}\cline{3-8}\cline{3-8}
 
	&&&\textit\textbf{{1 pt. move; top choice}} &  1.77 & 29.91 & 40.62 & 164.38 \\\cline{4-8}
        &&\textit{\textbf{Computation}}&\textit\textbf{{1 pt. move; top 2 choices}} &  1.99 & 61.6 & 105.12 & 436.97 \\\cline{4-8} 
	&&\textit{\textbf{Time (s)}}&\textit\textbf{{1 pt. move/swap; top choice}} & 3.85 & 82.43 & 150.47 & 595.81 \\\cline{4-8}
	&&&\textit\textbf{{1 pt. move/swap; top 2 choices}} & 2.11 & 158.83 & 352.16 & 3588.55 \\\hline
 
	\end{tabular}
\label{tab:percentage_by_neighborhood}
\end{table*}

For all these variations, the heuristic was run three times, and the best-obtained solution in terms of the objective value was obtained. The minimum, maximum, mean, and median percentage improvement with respect to the initial solution over the 43 instances are reported in Table~\ref{tab:percentage_by_neighborhood}. Further, in this table, the minimum, maximum, mean, and median computation time over the 43 instances are reported. The percentage improvement and computation time for all the instances for $\alpha = \frac{1}{TSP}, \tau = 1$ for all variations in the heuristic are shown in Fig.~\ref{fig:local_search_parametric_comp11}.
From the table and figure,
\begin{itemize}
    \item For $\alpha = \frac{1}{TSP}$ and $\tau = 1$ or $2,$ a median improvement of around $1-2 \%$, and a mean improvement of around $3-4\%$ was observed for all variations of the heuristic. Further, a maximum improvement of around $30\%$ was observed.
    \item Increasing $\alpha$ to $\frac{2}{TSP}$ leads to doubling the mean, median, and maximum improvements.
    \item The computation time for the heuristic for most cases is less than 10 minutes.
    \item Changing the neighborhood from 1 pt. move to 1 pt. move and swap leads to
    \begin{itemize}
        \item A marginal change in the improvement in objective value. For example, for $\alpha = \frac{1}{TSP}, \tau = 1,$ and top choice of maximal vehicle, the mean improvement increased by $0.13 \%,$ the median improvement increased by $0.37\%,$ and the maximum improvement increased by $0.98\%$.
        \item The computation time was about three times the previous computation time.
    \end{itemize}
    \item Changing the number of maximal vehicles considered from one vehicle to two vehicles leads to
    \begin{itemize}
        \item A larger change in the improvement in objective value. For example, for $\alpha = \frac{1}{TSP}, \tau = 1,$ and using 1 pt. move, the mean improvement increased by $0.48 \%,$ the median improvement increased by $0.64\%,$ and the maximum improvement increased by $4.34\%$.
        \item The computation time was about twice the previous computation time.
    \end{itemize}
    \item The choice of 1 pt. move with top two vehicles yields the best tradeoff between the objective value and computation time, whereas the choice of 1 pt. move and swap, and top two vehicles yields the best objective value.
\end{itemize}

\begin{figure}[htb!]
    \centering
    \subfigure[Percent Deviation of our solution wrt. \textit{memetic solution}]{\includegraphics[width = 0.35\textwidth]{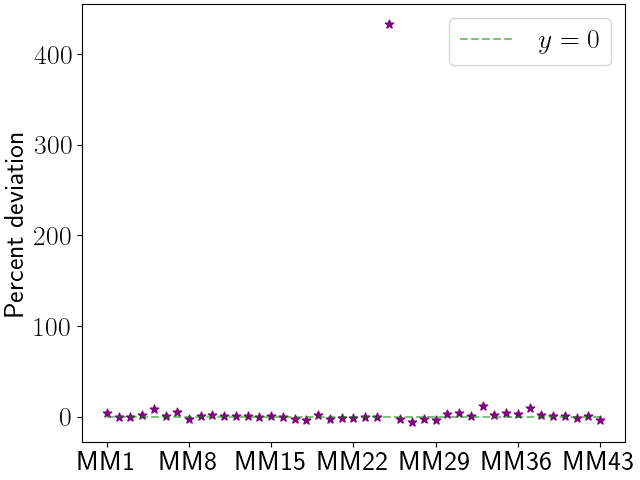}} \\
    \subfigure[Percent Deviation of our solution wrt. \textit{memetic solution} (without outlier)]{\includegraphics[width = 0.35\textwidth]{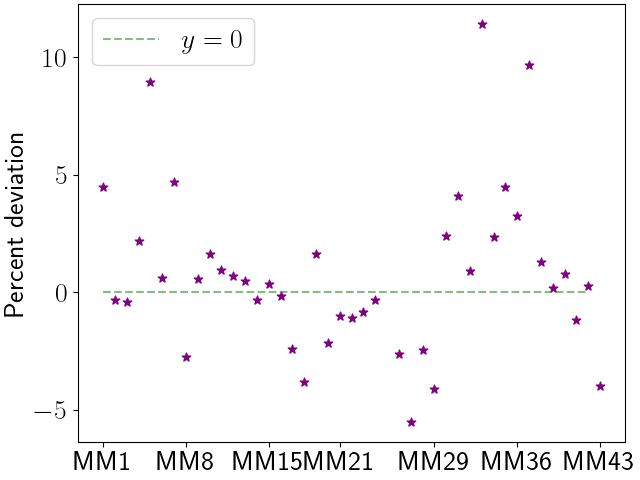}}
    \caption{Comparison of objective values from Min-Max tour solution vs. our algorithm using $\alpha=\frac1{TSP}$; $\tau=2$. The percentage difference has range$(-5.5\%,433.2\%)$ with a median of $0.35\%$. The mean and stdev. are $(10.9\%,65.3\%)$.
    }
    \label{fig:memetic_obj}
\end{figure}
\begin{figure}[htb!]
    \centering
        \subfigure[Memetic solution, objective value $=15.43$]{\includegraphics[width = 0.5\textwidth]{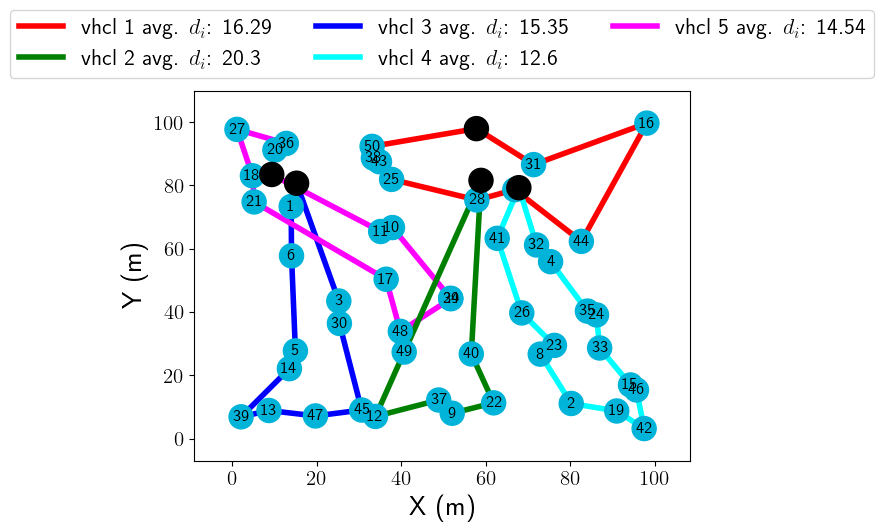}} \\
    \subfigure[Our solution, objective value $=15.69$]{\includegraphics[width = 0.5\textwidth]{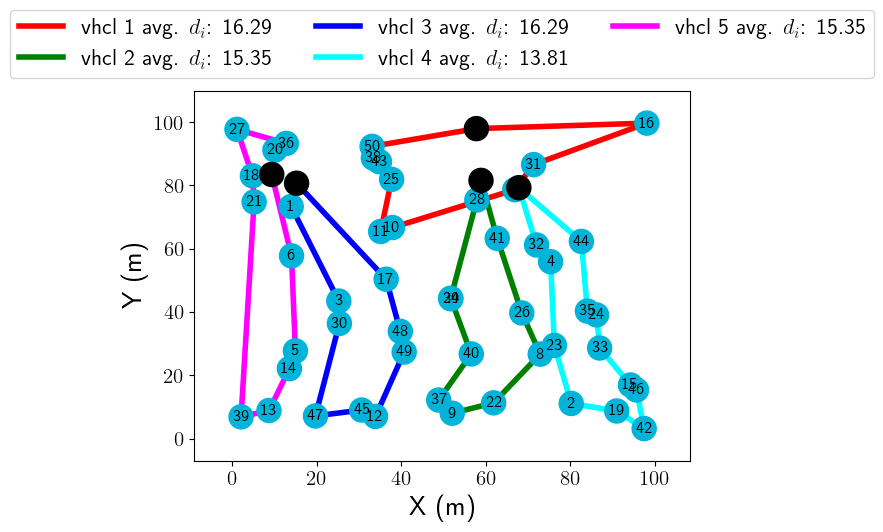}}
    \caption{Comparison of MM10 Min-Max tour solution vs. our algorithm solution for $\alpha=\frac1{TSP}=1.72\cdot 10^{-3}$; $\tau=2$ }
    \label{fig:memetic_paths}
\end{figure}

In addition, we compare our algorithm with results using tours obtained from the memetic algorithm \cite{memetic_algorithm}, the current best heuristic for the min-max multi-depot TSP. We utilize the partition of targets obtained from the memetic algorithm available at \cite{memetic_algorithm_solutions}. We compute the tour of each vehicle using LKH and the dwell time of targets covered by each vehicle using gradient descent. We choose to compare against the memetic algorithm since we expect the optimal solution to our problem to be similar to the min-max problem, since both problems require load balancing between vehicles. We expect this due to the results of Fig.~\ref{fig:targets}, where we see that the objective value gained by a vehicle roughly approaches a maximum before the discounting term takes over. 
For $\alpha = \frac{1}{TSP}$ and $\tau = 2$, 
the objective values obtained using each method is shown in Fig.~\ref{fig:memetic_obj}. For our heuristic, we used the 1 pt. move and swap, and top two choices of vehicles. We observe that the percentage difference of our solution from the solution using the partition from the memetic algorithm had a mean of $+11\%$ with a standard deviation of $65\%$, and a median difference of $0.35\%$. A depiction of solutions obtained using the heuristic and the memetic algorithm is shown for one of the instances in Fig.~\ref{fig:memetic_paths}, wherein the heuristic yields an improved objective value. We can see that on average, our algorithm performs marginally better than the solution obtained using the memetic algorithm. Further, we can see from Fig.~\ref{fig:memetic_obj} that there is one instance wherein our algorithm performs much better (up to a $400\%$ improvement), showing the benefit of our proposed heuristic for our proposed problem.

\section{Conclusions \& Remarks}

In this study, a single and multi-UAV routing problem for enhancing the performance of a classifier-in-the-loop system was considered, wherein an operator provides points-of-interest (POIs) through an interface, and vehicles are required to collect information regarding the same. The information gained was mathematically formulated using Kullback-Leibler divergence, and was discounted to ensure all POIs are visited. We considered two variants of the same problem: single-vehicle, and multi-vehicle. The single-vehicle problem was shown to be solved to optimality by decoupling the vehicle routing problem and optimizing the dwell time, which is the time spent at each POI. Numerical results for the same were presented to show that instances with 100 to 229 targets could be solved within 30 seconds. For the multi-vehicle variant, a heuristic was proposed due to coupled routing and dwell-time computation. Extensive numerical results were presented over varied types of instances along with variations in the heuristic to show that most of the instances and variations could be solved in 10 minutes. Furthermore, the heuristic was benchmarked with results obtained from a state-of-the-art heuristic for a min-max traveling salesman problem. It was observed that our proposed heuristic yielded an improved solution by about $11\%$ on average. Hence, the proposed heuristic produces high-quality solutions for the considered problem.

\section*{Acknowledgments}
Swaroop Darbha gratefully acknowledges the support of 2023 ONR Summer Faculty Fellowship Program. \\
\noindent{\bf DISTRIBUTION STATEMENT A:} This work is approved for public release (\# IR-5517-23-31-U). Distribution is unlimited. 

\bibliographystyle{IEEEtran}
\bibliography{cite}

\vfill

\end{document}